\documentclass[a4paper,10pt,reqno]{amsart}
\usepackage{amsthm,amssymb,latexsym,amsmath}
\usepackage[latin1]{inputenc}
\usepackage{amsmath,amssymb,amsfonts,amsthm}
\newtheorem{teo}{Theorem}[section]

\newtheorem{corollary}{Corollary}[section]

\newtheorem{defi}{Definition}[section]

\numberwithin{equation}{section}
\usepackage{color}

\newcommand{\F}{\mathcal F}
\newcommand{\E}{\mbox{${\mathcal E}$}}

\newcommand{\C}{\mbox{${\mathbb C}$}}

\title[ Rank two nilpotent co-Higgs sheaves  on complex surfaces]{Rank two nilpotent co-Higgs sheaves on complex surfaces}
\date{\today}
\author[ M. Corr\^ea]{M. Corr\^ea  }
\address{Departamento de Matem\'atica - ICEX, Universidade Federal de
  Minas Gerais, UFMG}
\curraddr{Av. Ant\^onio Carlos 6627, 31270-901,
  Belo Horizonte-MG, Brasil.}
\email{mauricio@mat.ufmg.br}
\begin{document}
\begin{abstract}
Let $(\E, \phi)$ be a rank two co-Higgs vector  bundles  on a K\"{a}hler  compact  surface $X$  with 
 $\phi\in H^0(X,End(\E)\otimes T_X)$ nilpotent. If    $(\E, \phi)$  is semi-stable,  then
one of the following holds up to finite  \' etale cover:
 \begin{itemize}
 \item[$i)$]  $X$ is uniruled. 
 \item[$ii)$]   $X$ is a torus and  $(\E, \phi)$ is strictly semi-stable.
\item[$iii)$] $X$  is a  properly elliptic surface and  $(\E, \phi)$ is strictly semi-stable.
\end{itemize}
\end{abstract}
\dedicatory{\it\ \ \ \ \ \ \ \ \ \ \ \ \ \ \ \ \ \ \ \ \ \ \ \ \ \ \ \  \ \ \ \ \ \ Dedicated to Jose Seade, for his 60th birthday.}
\maketitle
		
\section{Introduction.}
A generalised complex structure on a real manifold  $X$ of dimension $2n$, as defined by Hitchin \cite{H1}, is a rank-$2n$ isotropic subbundle
$E^{0,1}\subset (T_X\oplus T^*_X)^{\mathbb{C}}$ such that
\begin{itemize}
\item[i)] $E^{0,1}\oplus \overline{ E^{0,1} }=  (T_X\oplus T^*_X)^{\mathbb{C}}$
\item[ii)] $C^{\infty}(E^{0,1})$ is closed under the  Courant bracket.
\end{itemize} 
On a manifold with a generalized complex structure M.  Gualtieri in \cite{Gualtieri1} defined the notion of a generalized holomorphic bundle. More precisely,  a generalized holomorphic bundle on a generalized complex manifold, is a vector bundle $\E$ with a differential operator $\overline{D} : C^\infty(\E ) \longrightarrow C^{\infty}(\E \otimes E^{0,1})$ such that for all  smooth function $f$ and all  section $s\in C^\infty(\E) $ the following holds 
\begin{itemize}
\item[i)]$\overline{D} (fs)=\overline{\partial}(fs)+f \overline{D}(s)$
\item[ii)] $\overline{D}^2=0$.
\end{itemize} 
In the case of an ordinary complex structure and $\overline{D} = \overline{\partial}+\phi$, for operators
 $$\overline{\partial} : C^\infty(\E ) \longrightarrow C^{\infty}(\E \otimes \overline{T^*_X})$$
and
 $$\phi: C^\infty(\E ) \longrightarrow C^{\infty}(\E \otimes T_X),$$
the vanishing $\overline{D}^2=0$ means   that  $ \overline{\partial}^2=0$ ,  $ \overline{\partial} \phi=0$ and   $\phi\wedge\phi=0$.
By a classical result of Malgrange the condition $ \overline{\partial}^2=0$ implies that $\E$  is a holomorphic vector bundle. On the other hand,  $ \overline{\partial} \phi=0$ implies that $\phi$ is a holomorphic global section $$\phi\in H^0(X,End(\E)\otimes T_X)$$ which satisfies an  integrability condition   $\phi\wedge\phi=0$.
A co-Higgs sheaf  on a complex manifold $X$ is a sheaf  $\E$ together with a section  $\phi\in H^0(X,End(\E)\otimes T_X)$(called a  \textit{Higgs fields}) for which $\phi\wedge\phi=0$.
General properties of co-Higgs bundles were studied  in \cite{H,Ra2}.
There is a motivation in physics for studying co-Higgs bundles, see \cite{GW}, \cite{KLi} and \cite{Z}.

There are no stable co-Higgs bundles with nonzero Higgs field on curves $C$ of genus $g>1$. (When $g = 1$, a co-Higgs bundle is the same thing as a Higgs bundle in the usual sense.) In fact,   contracting with a holomorphic differential gives a non-trivial endomorphism of $\E$ commuting with $\phi$ which is impossible in the stable case \cite{H,Ra3} .
S. Rayan showed in \cite{Ra}  the non-existence  of  stable co-Higgs bundles with non trivial Higgs field on K3 and general-type surfaces. In this note we prove the following result.

 \begin{teo}\label{teo}
 Let $(\E, \phi)$ be a rank two co-Higgs vector  bundles  on a K\"{a}hler  compact  surface $X$  with 
 $\phi\in H^0(X,End(\E)\otimes T_X)$ nilpotent. If    $(\E, \phi)$  is semi-stable,  then
one of the following holds up to finite  \' etale cover:
 \begin{itemize}
 \item[$i)$]  $X$ is uniruled.
 \item[$ii)$]   $X$ is a torus and  $(\E, \phi)$ is strictly semi-stable.
\item[$iii)$] $X$  is a  properly elliptic surface and  $(\E, \phi)$ is strictly semi-stable.
\end{itemize}
 \end{teo}
 
 It follows direct  of proof of Theorem \ref{teo}  that  the we can consider a  more general classes of singular projective surfaces .
 
  \begin{teo}\label{teosing}
 Let $(\E, \phi)$ be a rank  two co-Higgs   torsion-free sheaf   on a normal projective  surface $X$  with 
 $\phi\in H^0(X,End(\E)\otimes T_X)$ nilpotent. If    $(\E, \phi)$  is stable,  then $X$ is uniruled.
 \end{teo}
 
 Finally, in this work we consider a relation  between 
co-Higgs bundles and Poisson geometry on  $\mathbb{P}^1$-bundles.  In  \cite{P} 
Polishchuk associated to  each  rank-$2$ co-Higgs bundle  $(\E, \phi)$  a Poisson structure on its  projectivized bundle $\mathbb{P}(\E)$. This relation was 
explained by Rayan in  \cite{Ra2} as follows:

 let $Y:= \mathbb{P}(\E)$  and consider the natural projection $\pi:Y \rightarrow X$. The exact sequence
$$
0\longrightarrow T_{X|Y}\longrightarrow T_Y\longrightarrow \pi^*T_X\longrightarrow 0
$$
 implies that $T_{X|Y}\otimes  \pi^*T_X \subset \bigwedge^2 T_Y$. Since $T_{X|Y}=\mathrm{Aut}(\mathbb{P}(\E))=\mathrm{Aut}(\E)/\mathbb{C}^*$ we get that
 $$\pi_*(T_{X|Y}\otimes  \pi^*T_X)=\pi_*T_{X|Y}\otimes T_X=End_0(\E)\otimes T_X,$$ where $End_0(\E)$ denotes the trace-free   endomorphisms   of $\E$. Therefore, we can  associate a trace-free co-Higgs fields  $\phi\in H^0(X,End(\E)\otimes T_X)$ a bi-vector  $\pi^*\phi\in H^0(X,T_{X|Y}\otimes  \pi^*T_X) \subset H^0(X, \bigwedge^2 T_Y)$ on $\mathbb{P}(\E)$. The co-Higgs condition $\phi\wedge\phi=0$ implies that  bi-vector  $\pi^*\phi$ is integrable, see the introduction of \cite{Ra2}. The codimension one foliation on $ \mathbb{P}(\E)$ is  the called    \textit{foliation by symplectic leaves } induced by Poisson struture . 
 
 We get an interisting consequence  of the proof Theorem \ref{teosing}. 
 
 \begin{corollary}\label{cor}
If     $(\E, \phi)$  is locally free, stable and nilpotent , then the closure of the all leaves of the foliation by symplectic leaves on   $\mathbb{P}(\E)$ 
 are   rational surfaces.
 \end{corollary}

\section{Semi-stable co-Higgs sheaves} 
 
\begin{defi}
A co-Higgs sheaf  on a complex manifold $X$ is a sheaf  $\E$ together with a section  $\phi\in H^0(X,End(\E)\otimes T_X)$(called a  Higgs fields) for which $\phi\wedge\phi=0$.
 \end{defi}
 
Denote by $End_0(\E) := ker(tr : End(\E) \longrightarrow \mathcal{O}_X)$  the trace-free part of the endomorphism bundle of $\E$.
  Since
 $$
 End(\E)= End_0(\E)\oplus \mathcal{O}_X
 $$
 we have that $End(\E)\otimes T_X=(End_0(\E)\otimes T_X)\oplus T_X$. Thus, the Higgs field $\phi\in H^0(X,End(\E)\otimes T_X)$ can be decomposed as $\phi=\phi_1+\phi_2$, where $\phi_1$ is the trace-free part and $\phi_2$ is a global vector field on $X$.  In particular, if the surface $X$ has no global holomorphic vector fields, then every Higgs field  
is trace-free.
\begin{defi}
Let $(X,\omega)$ be a polarized  K\"{a}hler  compact manifold. We say that $(\E,\phi)$ is semi-stable if
$$
\frac{c_1(\F)\cdot [\omega]}{rank(\F)}\leq \frac{c_1(\E)\cdot [\omega]}{rank(\E)}
$$
for all coherent subsheaves $0\neq \F\subsetneq \E$ satisfying $\Phi(\F)\subseteq \F\otimes T_X$, and stable  if the inequality  is strict for all such $\F$. We say that $(\E,\phi)$ is strictly semi-stable if  $(\E,\phi)$ is  semi-stable but non-stable.
\end{defi}
\section{Holomorphic foliations}
\begin{defi}\label{fol}
Let $X$ be a connected  complex manifold. A one-dimensional holomorphic foliation is
given by the following data:
\begin{itemize}
  \item[$i)$] an open covering $\mathcal{U}=\{U_{\alpha}\}$ of $X$;
  \item [$ii)$] for each $U_{\alpha}$ an holomorphic vector field $\zeta_\alpha$ ;
  \item [$iii)$]for every non-empty intersection, $U_{\alpha}\cap U_{\beta} \neq \emptyset $, a
        holomorphic function $$f_{\alpha\beta} \in \mathcal{O}_X^*(U_\alpha\cap U_\beta);$$
\end{itemize}
such that $\zeta_\alpha = f_{\alpha\beta}\zeta_\beta$ in $U_\alpha\cap U_\beta$ and $f_{\alpha\beta}f_{\beta\gamma} = f_{\alpha\gamma}$ in $U_\alpha\cap U_\beta\cap U_\gamma$.
\end{defi}
We denote by $K_{\F}$ the line bundle defined by the cocycle $\{f_{\alpha\beta}\}\in \mathrm{H}^1(X, \mathcal{O}^*)$. Thus, a one-dimensional holomorphic  foliation $\F$ on $X$ induces  a global holomorphic section $\zeta_{\F}\in \mathrm{H}^0(X,T_X\otimes K_{\F})$. The line bundle $K_{\F}$ is called the  \emph{canonical bundle} of $\F$.
Two sections  $\zeta_{\F}$ and $\eta_{\F}$ of $ \mathrm{H}^0(X,T_X\otimes K_{\F})$ are equivalent, if
there exists  a never vanishing holomorphic function $\varphi\in \mathrm{H}^0(X,\mathcal{O}^*) $,
such that $\zeta_{\F}= \varphi \cdot \eta_{\F}$.
It is clear that $\zeta_{\F}$ and $\eta_{\F}$ define the same foliation. Thus, a holomorphic foliation $\F$ on $X$  is an equivalence
of sections of $H^0(X, T_X\otimes K_{\F})$.
\section{Examples}
\subsection{Canonical example of split  co-Higgs bundles } 
Here we will give an example which naturally 
generalizes the canonical example given by Rayan on \cite[Chapter 6]{Ra2}.
Let $(X, \omega)$ be a polarized   K\"{a}hler  compact manifold. 
Suppose that there exists a global section $\zeta\in  
H^0(X , Hom(N,TX\otimes L))\simeq H^0(X ,TX\otimes L \otimes N^*)$. Now, consider the vector bundle
$$\E =L \oplus N.$$ 
Define the following  
co-Higgs fields 
$$\phi: L \oplus N \longrightarrow (T_X \otimes L)\oplus  (T_X \otimes N)\in H^0(X,End(\E)\otimes T_X)$$
by $\phi(s,t)=(\zeta(t),0)$. Since  $\phi\circ \phi=0 \in  H^0 (End(\E)\otimes T_X\otimes T_X) $ we get that  $\phi\wedge\phi=0$. Moreover, observe that the kernel of $\phi$  is the $\phi$-invariant line bundle $L$.
On the other hand, the line bundle  $L$ is destabilising only when 
$$ 
[2c_1(L)- c_1(\E)]\cdot [\omega]=[c_1(L)-c_1(N)]\cdot [\omega]> 0.
$$
That is, if 
$$ 
[c_1(L)]\cdot [\omega]> [c_1(N)]\cdot [\omega].
$$
\subsection{Co-Higgs bundles on ruled surfaces} 
Let $C$ be a curve  of genus $ g >1$. Now, consider the ruled surface $X:=\mathbb{P}(K_C\oplus \mathcal{O}_C)$ and $$\pi:\mathbb{P}(K_C\oplus \mathcal{O}_C)\longrightarrow C$$ the natural projection.  Consider  a Poisson structure
on $X$ given by a bivector $\sigma\in H^0(X, \wedge^2T_X)$. 
Let $(\E,\phi)$ be a nilpotent  Higgs bundle on $C$. S. Rayan showed   in \cite{Ra2}  that  $(\pi^*\E,\sigma(\pi^*\phi))$ is a stable co-Higgs bundle on $X$. 
\subsection{Co-Higgs orbibundles  on weighted projective spaces} 
Let $w_{0},w_1,w_{2}$ be positive integers, set $|w|:= w_0+ w_1 + w_2$. Assume that $w_{0},w_1,w_{2}$ are relatively prime.
Define an action of $\C^\ast$ in $\C^{3} \setminus \{0\}$ by
\begin{equation}
\begin{array}{ccc}
\C^\ast \times (\C^{3} \setminus \{0\} )& \longrightarrow &( \C^{3} \setminus \{0\}) \\
\lambda . (z_0, z_1, z_2) & \longmapsto &  (\lambda^{w_0} z_0, \lambda^{w_1} z_1, \lambda^{w_2} z_2)\\
\end{array}
\end{equation}
and consider the weighted projective plane
$$\mathbb{P}(w_{0},w_1,w_{2}) := (\C^{3} \setminus \{0\} )/ \sim$$ 
induced by the action above. We will denote this space  by  $\mathbb{P}(\omega)$.
On  $\mathbb{P}(\omega)$ we have an Euler sequence
$$
0\longrightarrow
\mathcal{O}_{\mathbb{P}(\omega)} \stackrel{\varsigma}{\longrightarrow}\bigoplus_{i=0}^{2
}\mathcal{O}_{\mathbb{P}(\omega)}(\omega_i) \longrightarrow
T\mathbb{P}(\omega) \longrightarrow 0,
$$
where $\mathcal{O}_{\mathbb{P}(\omega)}$ is the trivial line orbibundle  and  $T\mathbb{P}(\omega) = \mathrm{Hom} (\Omega_{\mathbb{P}(\omega)}^ 1,  \mathcal{O}_{\mathbb{P}(\omega)} )$ is the tangent  orbibundle of
$\mathbb{P}(\omega)$. The map $\varsigma$ is given explicitly by $
 \varsigma(1)=(\omega_0z_0,\omega_1z_1,\omega_2z_2).$ 
Now,  let $(\E,\phi)$ be a co-Higgs  orbibundle  on $\mathbb{P}(\omega)$.
Tensoring the Euler sequence by $End(\E)$, we obtain
$$ 0 \longrightarrow
End(\E) \longrightarrow
\bigoplus\limits_{i=0}^2 End(\E) ( \omega_i )
\longrightarrow  End(\E)  \otimes T \mathbb{P}(\omega)
 \longrightarrow 0.
$$
Thus, the co-Higgs fields $\phi$ can be represented, in homogeneous coordinates,  by
$$
\phi=\phi_0\otimes\frac{\partial}{\partial z_0}+\phi_1\otimes\frac{\partial}{\partial z_1}+\phi_2\otimes\frac{\partial}{\partial z_2} ,
$$
where $\phi_i\in H^0(\mathbb{P}(\omega), End(\E)(w_i))$, for all $i=0,1,2$,
and  $
\phi+\theta \otimes R_\omega $ define the same co-Higgs field as $\phi$, where $R_\omega$ is the adapted
radial vector field $$R_\omega = \omega_0 z_0 \frac{\partial}{\partial z_0}
+\omega_1 z_1 \frac{\partial}{\partial z_1}+   + \omega_2 z_2 \frac{\partial}{\partial z_2},$$ with $\theta$  a endomorphism of $\E$.
 Suppose that  
$$\E = \mathcal{O}(m_1) \oplus\mathcal{O}(m_2)$$ 
and that there exists a stable $\phi$ for $\E$ such that $m_1 \geq m_2$.
Then  $$|m_1- m_2| \leq \max\limits_{0 \leq i\neq j \leq 2}\{\omega_{i}+\omega_{j} \}.$$ 
In fact, this is a consequence of Bott's Formulae for weighted projective spaces. It follows from 
(see \cite{D})  that
\begin{center}$
\mathrm{H}^0(\mathbb{P}(\omega),T\mathbb{P}(\omega)\otimes\mathcal{O}_{\omega}(k))\simeq
\mathrm{H}^0(\mathbb{P}(\omega),
\Omega^{1}_{\mathbb{P}(\omega)}(\sum_{i=0}^{2} \omega_{i}+k))\neq
\emptyset$
\end{center}
if and only if $k>-\max\limits_{0 \leq i\neq j \leq 2}\{\omega_{i}+\omega_{j} \}$.
This generalize the example given by S. Rayan in \cite{Ra}. 
\subsection{Co-Higgs bundles on two dimensional complex tori} 
Let $X$ be a two dimensional complex torus and a co-Higgs bundle $\phi\in H^0(X,End(\E)\otimes T_X)$. Then  $\phi$ is equivalente to a pair of  commutative endomorphism  of $\E$. In fact,
since the tangent bundle $T_X$ is holomorphically trivial, we can take  a   trivialization  by choosing two  linearly independent global vector fields 
$v_1,v_2\in H^0(X,T_X)$. Then, we can write
$$
\phi= \phi_1 \otimes v_1 +\phi_2 \otimes v_2.
$$  
The condition $\phi\wedge\phi=0$ implies that 
$$
\phi_1\circ\phi_2=\phi_2\circ\phi_1.
$$
We have a  canonical nilpotent co-Higgs bundle $(\E,\phi)$, where  $\E=T_X=\mathcal{O}\oplus \mathcal{O}$ and 
$$\left({\begin{array}{*{20}c}
   0 & v  \\
   0 & 0  \\
 \end{array} } \right),
$$
where $v$ is a global vector field on $X.$
 
 \section{Proof of Theorem}

By using that  the Higgs field $\phi\in H^0(X,End(\E)\otimes T_X)$ is nilpotent we have  that  $Ker(\phi)=:L$ is a well defined line bundle on $X$ . Thus , we have a exact sequence 
  $$
  0\to L\longrightarrow  \E \longrightarrow \mathcal{I}_Z\otimes N \longrightarrow 0,
  $$
  where the nilpotent Higgs field $\phi$ factors as
  \begin{equation}
\begin{array}{ccc}
\E& \longrightarrow &\E \otimes T_X \\
\\
\downarrow &  &  \uparrow\\
\\
\mathcal{I}_Z\otimes N& \longrightarrow  &L\otimes T_X.
\end{array}
\end{equation}

The morphism $ \mathcal{I}_Z\otimes N\to L\otimes T_X $  induces a holomorphic foliation on $X$ which  induces 
 a global section  $\zeta_{\phi} \in H^0(X, T_X\otimes L\otimes N^*)$.   
Since $\det(\E)=L\otimes N$ we conclude that $L\otimes N^*=L^2\otimes \det(\E^*) $. Then
$$\zeta_{\phi}\in H^0(X, T_X\otimes L^2\otimes \det(\E^*) ).$$  
Let $K:=L^2\otimes \det(\E^*)$  the canonical bundle of the foliation $\F$ associated to the co-Higgs fields $\phi$.
If $\E$ is semi-stable then 
$$[c_1(K)]\cdot [\omega]=[2c_1(L)- c_1(\E)]\cdot [\omega]\leq 0
$$
for some K\"{a}hler  class $\omega$ .  If $K\cdot [\omega]<0$, it follows from \cite{Lamari} that $K$ is not pseudo-efective \cite{De}.
It follows from Brunella's theorem \cite{Br} that $X$ is uniruled. 
 
Now, suppose $X $ is a normal projective surface and  that $K\cdot H=0$, for some $H$ ample
 By Hodge index theorem we have that $K^2\cdot H^2\leq(K\cdot H)^2=0$, then $K^2\leq0$. 
Suppose that  $K^2<0$. We have that  $D=H+\epsilon  K$ is a $\mathbb{Q}$-divisor ample for $0<\epsilon<<1$, see \cite[proposition 1.3.6]{La}. Thus, we have that  
$$K\cdot D=K\cdot H+\epsilon^2K^2=\epsilon^2 K^2<0.$$ 
By  Bogomolov-McQuillan-Miyaoka's theorem  \cite{BMc} we conclude that  $X$ is uniruled.  
If $K^2=0$, then $K$ is numerically trivial. This fact  is well known, but for convenience of the reader we give a proof.  Suppose that there exists $C\subset X$ such that $K\cdot C>0$. Now, 
Consider the divisor $B=(H^2)C-(H\cdot C)H$. Then $B\cdot H=0$ and $K\cdot B=(H^2)K\cdot C<0.$ Define $F=mK+B$ for $0<m<<1$. Therefore 
$F\cdot H=0$ and $F^2>0$. This is a contradiction by the  Hodge index Theorem.
In this case $\E$ is strictly semi-stable.
  
Now, we apply the classification, up to finite \' etale cover,  of holomorphic foliations on projective surfaces with canonical bundle numerically trivial \cite{Peternel}, \cite{MQ2}, \cite{Br2}. 
Therefore, up to finite \' etale cover, either:
 \begin{itemize}
 \item[$i)$]  $X$ is uniruled; 
\item[$ii)$] $X$ is a torus;
\item[$iii)$] $k(X)=1$ and $X=B\times C$ with $g(B)\geq 2$,  $C$ is  elliptic. That is, $X$ is a  sesquielliptic surface.
\end{itemize}  
  
If $X$ is K\"{a}hler  and non-algebraic it follows from \cite{Br2} that,   up to finite \' etale cover, either 
$X$ has a unique elliptic fibration or $X$ is a torus.
 \section{Proof of Corollary \ref{cor}}
Since $(\E,\phi)$ is nilpotent and stable the co-Higgs fields induces a   foliation $\F$ by rational curves on $X$. 
Now, consider the projective bundle  $\pi:\mathbb{P}(\E) \rightarrow X$. Then the foliation by symplectic leaves $\mathcal{G}$ on   $\mathbb{P}(\E)$ 
is the pull-back of $\F$ by $\pi$. In particular, a  closure of the a leaf of the foliation by symplectic leaves $\mathcal{G}$ is of type  $\pi^{-1}(f(\mathbb{P}^1))$, where 
$f:\mathbb{P}^1\rightarrow X$ is the uniformization of a rational   leaf of $\F$. Clearly  $\pi^{-1}(f(\mathbb{P}^1))$ is a rational surface.
\\
\\
\\
\noindent{\footnotesize \textsc{Acknowlegments.} We are grateful to
Arturo Fernandez-Perez, Renato Martins and Marcos Jardim  for pointing out corrections.  We are grateful to Henrique Bursztyn for interesting conversations about Poisson Geometry.}

\enddocument